\def\Bbb#1{{\fam\msbfam\relax#1}}
\font\fivemsb=msbm5
\font\sevenmsb=msbm7
\font\tenmsb=msbm10
\newfam\msbfam
\textfont\msbfam=\tenmsb
\scriptfont\msbfam\sevenmsb
\scriptscriptfont\msbfam\fivemsb

\def\Bbb{\rm \bold}

\def\real{{\Bbb R}}
\def\complex{{\Bbb C}}

\def\vf{\varphi}

\def\wh{\widehat}
\def\a{\alpha}
\def\b{\beta}

\def\vF{\varPhi}

\def\f{\phi}

\def\m{\mu}

\def\r{\rho}
\def\s{\sigma}

\def\x{\xi}

\def\ce{{\cal E}}

\def\la{\langle}
\def\ra{\rangle}
\def\lla{\langle\!\langle}
\def\rra{\rangle\!\rangle}
\def\wh{\widehat}

\documentstyle[11pt,amssymb,amscd,righttag]{amsart}
\theoremstyle{plain}
\newtheorem{theorem}{Theorem}[section]

\newtheorem{proposition}[theorem]{Proposition}

\theoremstyle{definition}
\newtheorem{definition}[theorem]{Definition}
\newtheorem{example}[theorem]{Example}

\numberwithin{equation}{section}
\renewcommand{\rm}{\normalshape}
\begin{document}
\begin{flushright}
{\tt In: RIMS Kokyuroku {\bf 1139} (2000) pp19--29}
\end{flushright}
\bigskip

%\begin{flushright}
%{\tt typed Nov 21, 1999}
%\end{flushright}
%\medskip
%Topmatter
\title%[]%
{A Note on General Setting of White Noise Triple and 
Positive Generalized Functions}
%\thanks{The author is grateful for financial supports
%from the Daiko Foundation and the Kamiyama Foundation.}
\maketitle
%\author
\begin{center}
{\rm Nobuhiro ASAI\footnote{Current address: 
International Institute for Advanced Studies,
9-3 Kizu, Kyoto, 619-0225 Japan. 
E-mail address: asai@@iias.or.jp; nobuhiro.asai@@nifty.com}\\
Graduate School of Mathematics\\
Nagoya University\\
Nagoya, 464-8602 Japan \\
d97002h@@math.nagoya-u.ac.jp}
\end{center}
%\address{Nobuhiro Asai: Graduate School of Mathematics\\
%  Nagoya University \\ Nagoya 464-8602 \\ JAPAN}
\bigskip
\medskip
%\maketitle
%%%%%%%%%%%%%%%%%%%%%%%%%%%%%%%%%%%
% Introduction
%%%%%%%%%%%%%%%%%%%%%%%%%%%%%%%%%%%
\section{Introduction}
Let $\ce^{*}$ be the space of tempered distributions
and $\m$ be the standard Gaussian measure on $\ce^{*}$.
Being motivated by the distribution theory 
on infinite dimensional space
by Cochran, Kuo and Sengupta (CKS) \cite{cks}, 
Asai, Kubo and Kuo (AKK) have recently determined 
the best possible class $C_{+,{1\over 2},1}^{(2)}$
of functions $u$ to constract white noise triple,
\begin{equation*}
	[\ce]_u\subset L^2(\ce^{*},\m)
	\subset [\ce]^{*}_u,
\end{equation*}
and to characterize 
white noise test function space $[\ce]_u$ and generalized 
function space $[\ce]_u^{*}$ 
%from the viewpoint of the growth order
%of holomorphic functions (S-transform) 
in the series of papers 
\cite{akk1},\cite{akk2},\cite{akk3},\cite{akk4},\cite{akk5}.
The notion of Legendre transformation
plays important roles to examine relationships between 
the growth order of holomorphic functions
(S-transform) and the CKS-space of 
white noise test and generalized functions.  
It is well-known that a positive generalized 
function is induced by a Hida measure $\nu$ 
(generalized measure).
A Hida measure can be characterized 
by integrability conditions on a function 
inducing the above triple (\cite{akk5}). 
See also \cite{kuo99-1},\cite{kuo99-2},\cite{ob99} for 
an overview of other recent developments 
in white noise analysis. 

This short note is organized as follows.
In Section 2, we give a short summary of white noise analysis
including AKK's results.
A certain class of convex functions will be introduced to 
make use of Legendre transformation and dual functions
for our purposes.
In Section 3, we restate the characterization theorems 
of the spaces of white noise test and generalized functions
given in \cite{akk3},\cite{akk5}.
In Section 4, we give a quick review of 
the basic facts on the theory of positive generalized functions
\cite{ks76},\cite{kuo96},\cite{yokoi90}.
Finally, we discuss the characterization of a Hida measure 
(Theorem \ref{thm:mainthm}). 
In this connection, 
we present the grey noise and the Poisson noise 
measures as typical examples inducing 
positive generalized functions
(Examples \ref{example-grey} and 
\ref{example-poisson}, respectively). 
Moreover, we mention breifly the relationship
between $[\ce]_u$ and $L^s$-space on $(\ce^{*},\nu)$
(Proposition \ref{prop:L^p}).

%%%%%%%%%%%%%%%%%%%%%%%%%%%%%%%%%%%
% Preliminaries(Legendre and Triple)
%%%%%%%%%%%%%%%%%%%%%%%%%%%%%%%%%%%
\section{Preliminaries}
%\subsection{White Noise Space} 
Let us start with taking a special choice 
of a Gel'fand triple:
\begin{equation*}
	\ce={\cal S}(\real)\subset
	\ce_0=L^2(\real, dt)\subset
	\ce^{*}={\cal S}^{*}(\real)
\end{equation*}
just for convenience where 
${\cal S}$ is the Schwartz space of rapidly decreasing 
functions and ${\cal S}^{*}$ is the space of tempered
distributions.
Consult excellent books
\cite{kuo96},\cite{ob}
for more general constraction.
Let $A=1+t^2-d^2/dt^2$.  It is well-known that $A$
is a densely defined positive self-adjoint operator
on $\ce_0$.  So there exists an orthonormal
basis $\{e_j\}_{j=0}^{\infty}\subset\ce$ for $\ce_0$
satisfying $Ae_j=(2j+2)e_j$.  For each $p\geq 0$ 
we define $|f|_p=|A^pf|_0$ and let 
$\ce_p=\{f\in\ce_0 \ ; \ |f|_p<\infty, p\geq 0 \}$.
Note that $\ce_p$ is the completion of $\ce$
with respect to the norm $|\cdot|_p$.  Moreover, 
\begin{equation*}
	\r=\|A^{-1}\|_{op}={1\over 2}, 
	\quad \|i_{q,p}\|^2_{HS}
	=\sum_{j=0}^{\infty}(2j+2)^{-(q-p)}< \infty
\end{equation*}
for any $q>p\geq 0$.
Then the projective limit space $\ce$ of $\ce_p$ is
a nuclear space and the dual space
of $\ce$ is nothing but $\ce^{*}$.
Hence we have the following continuous inclusions:
\begin{equation*}
	\ce\subset\ce_p\subset\ce_0
	\subset\ce_p^{*}\subset\ce^{*},\quad p\geq 0.
\end{equation*}
where the norm on $\ce^{*}_p$ is given by 
\begin{equation*}
	|f|_{-p}=|A^{-p}f|_0.
\end{equation*}
Troughout this paper, we denote 
the complexification of a real space $X$
by $X_c$.
Let $\m$ be the standard Gaussian
measure on $\ce^{*}$ given by 
\begin{equation*}
	e^{-{1\over 2}|\xi|^2_0}
	=\int_{\ce^{*}}e^{i\langle x,\xi\rangle}\mu (dx),
	\quad\quad \xi\in\ce^{*}.
\end{equation*}
The probability space $(\ce^{*},\m)$ is called the 
{\it white noise space} or {\it Gaussian  space}.
Let $(L^{2})=L^2(\ce^{*},\m)$ denote the Hilbert space of
$\m$-square integrable functions on $\ce^{*}$. By the
Wiener-It\^o theorem each $\vf$ in $(L^{2})$ can be uniquely
expressed as
\begin{equation} \label{eq:1-1}
	\vf = \sum_{n=0}^{\infty} I_{n}(f_{n}) = \sum_{n=0}^{\infty}
	\la :\!\cdot^{\otimes n}\!:, f_{n}\ra, \qquad
	f_{n}\in \ce_{0,c}^{\wh\otimes n},
\end{equation}
and the $(L^{2})$-norm $\|\vf\|_{0}$ of $\vf$ is given by
\begin{equation}
	\|\vf\|_{0} 
	= \left(\sum_{n=0}^{\infty} n!|f_{n}|_{0}^{2}\right)^{1/2}.   
	\notag
\end{equation}

\smallskip
We now briefly recall notions and results discussed by
Asai et al. \cite{akk2},\cite{akk3}.
\smallskip
\begin{definition}
Let $u$ be a positive continuous function on $[0, \infty)$.
%\begin{itemize}
%\item[(a)] The function $u$ is called {\em log-convex} if
%$\log u$ is convex on $[0, \infty)$.
%\item[(b)] The function $u$ is called {\em (log, exp)-convex}
%if $\log u(e^{x})$ is convex on $\spr$.
%\item[(c)] 
The function $u$ is called
{\em (log, $x^{2}$)-convex} if $\log u(x^{2})$ is convex on
$[0, \infty)$. 
%Here $k$ is a positive real number.
%\end{itemize}
\end{definition}  

\begin{definition}
\begin{itemize}
\item[(1)] Let $C_{+,\log}$ be the class of all positive
continuous functions on $[0,\infty)$ and 
\begin{equation}\label{eq:entire}
	\lim_{r \to\infty}{\log u(r) \over \log r}=\infty.
\end{equation}
\item[(2)] Let $C_{+,{1\over 2}}$ be the class of all positive
continuous functions on $[0,\infty)$ and 
\begin{equation}\label{eq:entire}
	\lim_{r \to\infty}{\log u(r)\over\sqrt{r}}=\infty.
\end{equation}
\item[(3)] Let $C_{+,{1\over 2},1}$
be the class of functions $u\in C_{+,{1\over 2}}$
such that
\begin{equation}\label{eq:embed}
	\lim_{r\to\infty}{\log u(r)\over r}<\infty.
\end{equation}
\item[(4)] Let $C^{(2)}_{+,{1\over 2}}$ \
( or $C^{(2)}_{+,{1\over 2},1}$)
be the class of all $(\log, x^2)$-convex functions
in $C_{+,{1\over 2}}$ \ (or $C_{+,{1\over 2},1}$),
respectively.
%and $C^{(\infty)}_{+,{1\over 2}}$ \ 
%( or $C^{(\infty)}_{+,{1\over 2},1}$) 
%be the class of all $(\log, \exp)$-convex functions
%in $C_{+,j}$ \ (or $C_{+,j,1}$)
%\item[(5)]Let $C^{(k)}_{+,\log}$ \
%be the class of all $(\log, x^k)$-convex functions
%in $C_{+,\log}$ 
%and $C^{(\infty)}_{+,\log}$ \ 
%be the class of all $(\log, \exp)$-convex functions
%in $C_{+,\log}$.
\end{itemize}
\end{definition}

\begin{definition}
The {\it Legendre transform} $\ell_u$
of a function
$u\in C_{+,\log}$ is defined by
\begin{equation}\label{eq:ell}
	\ell_u(t):=\inf_{r>0}{u(r)\over r^t},
\quad t\geq 0
\end{equation}
\end{definition}
%and its inverse transform is given by
%\begin{equation}\label{eq:inv-formula}
%	u(r)=\sup_{t\geq 0}\ell_u(t)r^t
%	\ \ \hbox{for} \ r>0.
%\end{equation}

\begin{definition}
The {\it dual Legendre function}
$u^{*}$ of a function $u\in C_{+,{1\over 2}}$ 
is given by
\begin{equation}\label{eq:dualfunction}
	u^{*}(r):=\sup_{s>0}{e^{2\sqrt{rs}}\over u(s)}
	\in C^{(2)}_{+,{1\over 2}}.
\end{equation}
\end{definition}

\begin{definition}
	Two positive sequences
	$\{a(n)\}_{n=0}^{\infty}$ 
	and 
	$\{b(n)\}_{n=0}^{\infty}$ 
	are called {\it equivalent} 
	(denoted by $a(n)\sim b(n)$) if 
	there exist positive constants 
	$K_1$, $K_2$, $c_1$, $c_2$ such that  
\begin{equation}
  	K_1c_1^na(n)\leq b(n) \leq K_2c_2^na(n)
  	\ \ \hbox{for any} \ n\in {\Bbb N}.
\end{equation} 
\end{definition}

\begin{definition}
	A positive sequence 
	$\{a(n)\}_{n=0}^{\infty}$ 
	is a dual sequence of 
	$\{b(n)\}_{n=0}^{\infty}$
	if $a(n)b(n)\sim (n!)^{-2}$ holds.
\end{definition}
According to Theorem 4.6 in \cite{akk3},
for a fucntion
$u\in C_{+,{1\over 2}}^{(2)}$
we have 
\begin{equation}
	\ell_{u^{*}}(n)
	\sim {1\over \ell_u(n)(n!)^2}.
\end{equation}
Thus, $\{\ell_{u^{*}}(n)\}_{n=0}^{\infty}$ 
is a dual sequence of 
$\{\ell_u(n)\}_{n=0}^{\infty}$.

\begin{definition}
	Two positive functions 
	$u(r)$ and $v(r)$ on $[0,\infty)$ 
	are called {\it equivalent} 
	(denoted by $u(r)\approx v(r)$) if 
	there exist positive constants 
	$a_1$, $a_2$, $c_1$, $c_2$ such that  
\begin{equation}
  	c_1u(a_1r)\leq v(r) \leq c_2u(a_2r)
  	\ \ \hbox{for any} \ r\geq 0.
\end{equation}
\end{definition}

%The $(\log,x^2)$-convexity condition ensures one to one 
%correspondence between $u$ and $\ell_u(n)$ ($u$ and $u^{*}$).
The condition \eqref{eq:entire} is required 
for both $u$ and $u^{*}$ to be equivalent 
to entire functions, respectively.  
This requirement is essential 
for Theorem \ref{thm:test-chara} 
and Theorem \ref{thm:dis-chara}
which will be discussed later.

Next, we describe the spaces of test and generalized functions
on the space $\ce^{*}$ introduced by Cochran et al.~in a recent
paper \cite{cks}. Let $\{\a(n)\}_{n=0}^{\infty}$ be a
weight sequence satisfying the following 
{\it two conditions} \cite{akk2},\cite{akk4},\cite{cks}:
\begin{itemize}
\item[(A1)] 
$\a(0)=1$ 
and $\inf_{n\geq 0}\,\a(n)\s^n> 0$ for some $\s\geq 1$.
\item[(A2)] 
$\lim_{n\to\infty} 
\left({\a(n)\over n!}\right)^{1/n} =0$.
\end{itemize}

Let $\vf\in (L^{2})$ be represented 
as in Equation (\ref{eq:1-1}). 
For $p\geq 0$ and  
a given function $u\in C^{(2)}_{+,{1\over 2},1}$,
define
\begin{equation} \label{eq:a}
	\|\vf\|_{u,p} 
	= \left(\sum_{n=0}^{\infty}{1\over \ell_u(n)}
	|f_{n}|_{p}^{2} \right)^{1/2},
	\qquad f_n\in\ce^{\wh{\otimes}n}_{p,c}.
\end{equation}
where $|f_n|_p=|(A^{\otimes n})^pf_n|^2_0$.
Note that since 
$\ell_u(n)$ and $\ell_{u^{*}}(n)$ are 
dual sequences of one another, we choose 
$\a(n)=(n!\ell_u(n))^{-1}\sim n!\ell_{u^{*}}(n)$ 
as a weight sequence.  
It is easy to check 
by \eqref{eq:entire}, \eqref{eq:embed} and \eqref{eq:ell}
that the above weight sequence 
with assumptions, $u\in C_{+,{1\over 2},1}^{(2)}$ 
and $\inf\{u(r); r\geq 0\}=1$, satisfies
the conditions (A1) and (A2).  
Let 
	$[\ce_{p}]_{u} 
	= \{\vf\in (L^{2})\, ;\, \|\vf\|_{u,p}<\infty\}$. 
Define the space $[\ce]_{u}$ of test functions
to be the projective limit of the family 
$\{[\ce_{p}]_{u};\,p\geq 0\}$. 
Its dual space $[\ce]_{u}^{*}$ is the space of
generalized functions. By identifying $(L^{2})$ with its
dual we get the following continuous inclusion maps:
\begin{equation}
	[\ce]_{u} \subset [\ce_{p}]_{u} \subset (L^{2}) \subset
	[\ce_{p}]_{u}^{*} \subset [\ce]_{u}^{*},
	\quad p\geq 0. \notag
\end{equation}
Note that the condition \eqref{eq:embed} is needed in order to 
have the continuous inclusion 
$[\ce_p]_{u}\subset (L^2)$ for $p\geq 0$.
The canonical bilinear form on $[\ce]^*_{u}\times [\ce]_{u}$
is denoted by $\lla\cdot,\cdot\rra$.
For each $\varPhi\in [\ce_p]^{*}_{u}$ there exists a unique
$F_n\in (\ce^{\otimes n}_{p,c})^*_{symm}$ 
such that 
\begin{equation*}
  	\lla\varPhi,\varphi\rra=\sum_{n=0}^{\infty}n!\la F_n,f_n\ra
\end{equation*}
and
\begin{equation}  \label{eq:a-1}
	\|\vf\|_{u^{*},-p} 
	= \left(\sum_{n=0}^{\infty} 
	{1\over \ell_{u^{*}}(n)}
	|F_{n}|_{-p}^{2} \right)^{1/2}.
\end{equation}
The Gel'fand triple $[\ce]_{u} \subset (L^{2}) \subset
[\ce]_{u}^{*}$ is called the {\it Cochran-Kuo-Sengupta space} 
(CKS-space for short)
associated with a given weight function 
$u\in C_{+,{1\over 2},1}^{(2)}$
(see \cite{akk4},\cite{akk5}).
We remark that some similar results have been obtained
independently by Gannoun et al. \cite{ghor}.
\smallskip

%%%%%%%%%%%%%%%%%%%%%%%%%%%%%%%%%%
% Examples
%%%%%%%%%%%%%%%%%%%%%%%%%%%%%%%%%%
%\newpage
%\smallskip
\begin{example}\label{example2.8}
Let $0\leq\beta <1$.
If $v(r)=\exp((1+\beta)r^{{1\over 1+\beta}})$,
then $v\in C^{(2)}_{+,{1\over 2},1}$.
In addition, the dual function is 
$v^{*}(r)=\exp((1-\beta)r^{{1\over 1-\beta}})$.
By Stirling formula, we have for any $n\geq 0$
\begin{equation}
	{1\over (n!)^{1+\b}}
	\leq \ell_{v}(n)=\biggl({e\over n}\biggr)^{(1+\b)n}
	\leq \biggl({e2^{n/2}\over n!}\biggr)^{1+\b}
\end{equation}
That is, $\ell_{v}(n)\sim (n!)^{-(1+\b)}$.
On the other hand, similarly we get 
$\ell_{v^{*}}(n)\sim (n!)^{-(1-\b)}$.
Hence $\ell_{u_{1}}(n)$ and $\ell_{v^{*}}(n)$
are dual sequences of one another.
This example induces the {\it Hida-Kubo-Takenaka space} 
with $\beta=0$ \cite{hkps},\cite{kt1},\cite{kt2},\cite{ob}, 
$$
  	(\ce)_0\subset (L^2)\subset (\ce)^*_0,
$$
and the {\it Kondratiev-Streit space} 
\cite{ks92},\cite{ks93},\cite{kuo96},
$$
  	(\ce)_{\beta}\subset (L^2)\subset (\ce)^*_{\beta},
  	\ \ 0\leq \beta < 1.
$$
See \cite{houz}, \cite{kls96} if $\beta=1$. 
\end{example}

\begin{example}\label{example2.9}
Let 
$$
	u_k(r)={\exp_k(r)\over\exp_k(0)}
	=\sum_{n=0}^{\infty}{b_k(n)\over n!}r^n. 
$$
where 
\begin{equation*}
  	\exp_1(r):=\exp(r),
  	\ \ \exp_k(r):=\exp(\exp_{k-1}(r))
  	\ \ \hbox{for}\ k\geq 2
\end{equation*}
and $b_k(n)$ is {\it the k-th order Bell number} 
\cite{akk1},\cite{akk2},\cite{akk3},\cite{akk4},\cite{cks},\cite{kubo}. 
Then its dual Legendre function $u^{*}_2$ is 
equivalent to the function 
$$
	u^{*}_k(r)\approx\exp\biggl[2\sqrt{r\log_{k-1}\sqrt{r}}\biggr],
$$
where $\log_k(r)$ is given by 
\begin{equation*}
  	\log_1(r):=\log(\max\{e,r\}),
  	\ \ \log_{k}(r):=\log_1(\log_{k-1}(r))
  	\ \ \hbox{for}\ k\geq 2.
\end{equation*}
Then, $u^{*}_k\in C^{(2)}_{+,{1\over 2},1}$.
Details of $\ell_{u_{k}}(n)$
can be found in \cite{kubo},\cite{kks}.
The Gel'fand triple
$$
  	[\ce]_{u_k}
  	\subset (L^2)
  	\subset [\ce]^{*}_{u_k}
$$
is called the {\it Bell number space of order $k$}
\cite{cks} and 
$$
  	[\ce]_{u_k}\subset(\ce)_{\beta}
  	\subset (L^2)\subset (\ce)^*_{\beta}
  	\subset [\ce]^{*}_{u_k}.
$$
\end{example}

\smallskip
\noindent
{\bf Remark.} 
We point out here that 
\begin{equation*}
  	(\ce)_1 \subset [\ce]_{u}\subset (\ce)_0\subset (L^2)
  	\subset (\ce)^*_0 \subset [\ce]^*_{u} \subset (\ce)^*_1
\end{equation*}
holds for any $u\in C^{(2)}_{+,{1\over 2},1}$. \\
%%%%%%%%%%%%%%%%%%%%%%%%%%%%%%%%%%%%%%%%%%%%%%%%
% Characterization
%%%%%%%%%%%%%%%%%%%%%%%%%%%%%%%%%%%%%%%%%%%%%%%
\section{General Characterization Theorems}
The 
%renormalized exponential function
{\it exponential function} (or {\it coherent state}) 
$\f_{\x}(\cdot)$
%=:\!e^{\la\cdot, \x\ra}\!:$
is given by 
\begin{equation*}
%:\!e^{\la\cdot, \x\ra}\!: \
	\f_{\x}
	=e^{\la\cdot,\x\ra-{1\over 2}|\x|^2_0}
	=\sum_{n=0}^{\infty}{1\over n!}
	\la :\!\cdot^{\otimes n}\!:, \xi^{\otimes n}\ra,
	\qquad \xi\in\ce_c.
\end{equation*}
Since it is well-known that the exponential functions
$\{
%:\!e^{\la\cdot, \x\ra}\!: 
\f_{\x}
\ ; \ \xi\in\ce_{c}\}$
span a dense subspace of $[\ce]_{u}$, $\varPhi\in [\ce]^{*}_{u}$
is uniquely determined by its {\it S-transform}:
\begin{equation*}
	S\varPhi(\xi):=\lla\varPhi, 
	%:\!e^{\la\cdot, \x\ra}\!: 
	\f_{\x} \rra,
	\ \ \xi\in\ce_{c}.
\end{equation*}

%\noindent
\smallskip
\noindent
Now we are in a position to state
the characterization theorems of
test and generalized functions 
under very general assumptions.
We remark that these theorems were examined in  
\cite{kps},\cite{ps} for the Hida-Kubo-Takenaka
space case and 
\cite{ks92},\cite{ks93} for the 
Kondratiev-Streit space case.
%%%%%%%%%%%%%%%%%%%%%%%%%%%%%%%%%%%%%%%%%%%%%%%%%%
% Characterization of test functons
%%%%%%%%%%%%%%%%%%%%%%%%%%%%%%%%%%%%%%%%%%%%%%%%%%
\begin{theorem}
[\cite{akk3},\cite{akk5}]\label{thm:test-chara}
Let $u\in C^{(2)}_{+,{1\over 2},1}$
be increasing with $u(0)=1$.
A function $F: \ce_{c}\rightarrow \Bbb C$
is the S-transform of some $\varphi\in [\ce]_{u}$ if and only if \\
$(F1) z\mapsto F(z\xi+\eta)$ is entire holomorphic in
$z\in\Bbb C$ for any $\xi, \eta\in \ce_{c},$ \\
$(F2)'$ For any $a, p\geq 0$, there exists a constant $K>0$ such that
\begin{equation}
	|F(\xi)|\leq Ku(a|\xi|^2_{-p})^{{1\over 2}}
	\ \ \hbox{for any}\ \xi\in \ce_{c}.
\end{equation}
In addition, in that case,
\begin{equation}
	\|\varphi\|^2_{u,q}
	\leq K^2\bigl(1-ae^2\|i_{p,q}\|^2_{HS}\bigr)^{-1}
\end{equation}
for any $q<p$ satisfying $ae^2\|i_{p,q}\|^2_{HS}<1$.
\end{theorem}

%%%%%%%%%%%%%%%%%%%%%%%%%%%%%%%%%%%%%%%%%%%
% Characterization of Generalized functions
%%%%%%%%%%%%%%%%%%%%%%%%%%%%%%%%%%%%%%%%%%
\begin{theorem}[\cite{akk5}]\label{thm:dis-chara}
Suppose that $u\in C^{(2)}_{+,{1\over 2},1}$
and $\inf_{r>0}u(r)=1$.
A function $F: \ce_{c}\rightarrow \Bbb C$
is the S-transform of some 
$\vf\in [\ce]^*_{u}$ if and only if \\
$(F1) z\mapsto F(z\xi+\eta)$ is entire holomorphic in
$z\in\Bbb C$ for any $\xi, \eta\in \ce_{c},$ \\
$(F2)$ there exist nonnegative constants $K,a$ and $p$ such that
\begin{equation}
	|F(\xi)|\leq Ku^{*}(a|\xi|^2_{p})^{{1\over 2}}
	\ \ \hbox{for any}\ \xi\in \ce_{c}.
\end{equation}
In addition, in that case,
\begin{equation}
	\|\vf\|^2_{u^{*},-q}
	\leq K^2\bigl(1-ae^2\|i_{q,p}\|^2_{HS}\bigr)^{-1}
\end{equation}
for any $q>0$ satisfying $ae^2\|i_{q,p}\|^2_{HS}<1$.
\end{theorem}

\noindent
{\bf Remark.}
In \cite{hida75}, Hida introduced the
infinite dimensional analogue of Fourier transform, so-called   
{\it ${\cal T}$-transform}, given by 
\begin{equation}
	{\cal T}\vF(\x)=\lla\vF, e^{i\la\cdot,\x \ra}\rra,
	\quad \x\in\ce_c.
\end{equation}
It is well-known that 
the ${\cal T}\vF(z\x+\eta)$ is entire holomorphic
in $z\in \complex$ for any $\x,\eta\in\ce_{c}$.
In addition, there is a nice relationship 
between ${\cal T}$-transform and $S$-transform:
\begin{equation}\label{eq:ST-relation1}
	S\vF(\x)={\cal T}(-i\x)\exp [-{1\over 2}\la\x,\x\ra],
	\quad \x\in\ce_c.
\end{equation}
and 
\begin{equation}\label{eq:ST-relation2}
	{\cal T}\vF(\x)=S(i\x)\exp [-{1\over 2}\la\x,\x\ra],
	\quad \x\in\ce_c.
\end{equation}
Therefore, Theorems \ref{thm:test-chara} and \ref{thm:dis-chara}
remain valid even if the $S$-transform is replaced by 
the ${\cal T}$-transform.

%%%%%%%%%%%%%%%%%%%%%%%%%%%%%%%%%%%%%%%%%%%%%%%%%%%%%%%%%%%%%%
% Hida measure
%%%%%%%%%%%%%%%%%%%%%%%%%%%%%%%%%%%%%%%%%%%%%%%%%%%%%%%%%%%%%%
\section{Characterization of Positive Generalized Functions}

\begin{definition}
A measure $\nu$ on $\ce^{*}$ is called a {\it Hida measure}
if $[\ce]_{u}\subset L^1(\nu)$ and the mapping 
$\vf\mapsto \int_{\ce^{*}}\vf(x)\nu(dx)$ is coninuous on $[\ce]_u$.
\end{definition}

Suppose $\vF\in [\ce]^{*}_u$ is induced by a 
Hida measure $\nu$ on $\ce^{*}$.  Then 
\begin{equation*}
	\lla \vF, \vf\rra
	=\int_{\ce^{*}}\vf(x)\nu(dx), \ \vf\in [\ce]_u.
\end{equation*}
In paticular, take $\vf=e^{i<x,\x>}, \ \x\in\ce$.
Then, we have the ${\cal T}$-transform restricted to
$\ce$ by   
\begin{equation*}
	{\cal T}\vF(\x)|_{\ce}:=\lla \vF, e^{i<x,\x>}\rra
	=\int_{\ce^{*}}e^{i<x,\x>}\nu(dx), \quad \x\in\ce.
\end{equation*}
It is clearly seen that 
${\cal T}\vF(\x)|_{\ce}$ is equal to the characteristic 
function $C(\x)$ of $\nu$.
Thus in order to see  
the existance of a measure $\nu$, 
we only need to check that the function
${\cal T}\vF(\x)|_{\ce}$ satisfies the following 
conditions: \\
(1) ${\cal T}\vF(\x)|_{\ce}$ is continuous on $\ce$. \\
(2) ${\cal T}\vF(\x)|_{\ce}$ is positive definite on $\ce$.

\begin{definition}
A generalized function 
$\vF\in [\ce]_u^{*}$
is called {\it positive} if 
$\lla \vF, \vf\rra\geq 0$
for all nonnegative test functions 
$\vf\in[\ce]_u$. 
\end{definition}

Hence for any nonnegative test function $\vf\in [\ce]_u$,
\begin{equation*}
	\lla \vF, \vf\rra
	=\int_{\ce^{*}}\vf(x)\nu(dx)
	\geq 0.
\end{equation*}
Thus $\vF$ is a positive generalized function.
In the following, Theorem \ref{thm:3-equiv} says 
that all of positive generalized functions 
is generalized functions induced by Hida measures on $\ce^{*}$.
\begin{theorem}\label{thm:3-equiv}
Let $u\in C_{+,{1\over 2},1}^{(2)}$ and
$\vF\in [\ce]^{*}_{u}$.  Then the following 
statements are equivalent: 
\begin{itemize}
\item[(a)] $\vF$ is positive. 
\item[(b)] ${\cal T}\vF$ is positive definite on $\ce$. 
\item[(c)] $\vF$ is induced by a Hida measure.  That is,
there exists a finite measure $\nu$ on $\ce^{*}$
such that $[\ce]_u\subset L^1(\nu)$ and 
\begin{equation*}
	\lla \vF, \vf\rra =
	\int_{\ce^{*}}\vf(x)\nu(dx), \quad \vf\in [\ce]_u.
\end{equation*}
\end{itemize}
\end{theorem}

\noindent
{\bf Remark.}
For the Kondratiev-Streit space, 
the equivalence of $(a)\sim (c)$ has been discussed
in \cite{kuo96}.
The equivalence of $(a)$ and $(c)$ was examined originally in 
\cite{ks76} and \cite{yokoi90} only 
for the case of the Hida-Kubo-Takenaka space.

\smallskip
\smallskip
Next, our problem is how to characterize Hida measures on $\ce^{*}$.
For this purpose, based on Lee's idea \cite{lee} 
(see also \cite{akk3},\cite{akk5},\cite{kuo96},\cite{oue}), 
we shall define another norm as follows.
Let ${\cal A}_{u,p}$
be the space of all functions
$\varphi$ on $\ce^{*}_{c}$ 
satisfying the following conditions:\\
(A1) $\varphi$ is an analytic function on 
$\ce^{*}_{p,c}$. \\
(A2) There exists a nonnegative constant $C$ such that
\begin{equation*}
	|\varphi (x)|^2\leq Cu(|x|^2_{-p})
	\ \ \hbox{for any}\ 
	x\in \ce^{*}_{p,c}.
\end{equation*}
For $\varphi\in {\cal A}_{u,p}$, 
its norm is defined by
\begin{equation}\label{eq:Anorm}
	|\!|\!|\vf|\!|\!|_{{\cal A}_{u,p}}
	:=\sup_{x\in \ce^{*}_{p,c}}
	|\varphi (x)|u(|x|^2_{-p})^{-{1\over 2}}.
\end{equation}
for a function $u\in C_{+,\log}$.
Define the space ${\cal A}_{u}$ of
test functions on $\ce^{*}$ by
\begin{equation*}
	{\cal A}_{u}
	=\projlim_{p\to\infty}
	{\cal A}_{u,p}.
\end{equation*}

\begin{proposition}
Suppose that $u\in C_{+,{1\over 2},1}^{(2)}$.
Then $\{|\!|\!|\cdot|\!|\!|_{{\cal A}_{u,p}}; p\geq 1\}$
is equivalent to $\{\|\cdot\|_{u,p}; p\geq 1 \}$.
As a result, ${\cal A}_u=[\ce]_u$ as vector spaces
\end{proposition}
\noindent
{\bf Remark.}
This proposition implies that 
$$
	{\cal A}^*_{u}=[\ce]^*_{u}
$$
and they have the same inductive limit topology.
Moreover, for the constraction of spaces ${\cal A}_u$ 
and ${\cal A}^{*}_u$ the Wiener-It\^o
decomposition theorem is not used and 
a measure on $\ce^{*}$ is not refered at all.

\smallskip
We are in a position to state 
our characterization theorem of Hida measure
in terms of positivity and integrability condition.  
The proof of Theorem \ref{thm:mainthm} is 
based on simple applications of 
$S$-transform (${\cal T}$-transform, 
equivalently the Bochner-Minlos Theorem),
the dual functions given in \eqref{eq:dualfunction}
and some technical estimations.
See \cite{akk5} for details.
 
\begin{theorem}[\cite{akk5}]\label{thm:mainthm}
Suppose that $u\in C^{(2)}_{+,\frac{1}{2},1}$.
A measure $\nu$ 
on $\ce^{*}$ is a Hida measure
inducing a positive generalized function
$\varPhi_{\nu}\in [\ce]_u^{*}$
if and only if $\nu$ 
is supported in $\ce^{*}_{p}$
for some $p\geq 1$ and
\begin{equation*}
	\int_{\ce^{*}_{p}}
	u(|x|^2_{-p})^{\frac{1}{2}}
	\nu(dx)<\infty.
\end{equation*}
\end{theorem}

\noindent
{\bf Remark.}
See \cite{lee} for the Hida-Kubo-Takenaka space
and \cite{kuo96} for the Kondratiev-Streit space.

\begin{proposition}\label{prop:L^p}
Let $u\in C_{+,{1\over 2},1}^{(2)}$ and 
$\nu$ be a Hida measure on $\ce^{*}$
inducing a generalized function $\vF_{\nu}$ in $[\ce]^{*}_u$.
Then $[\ce]_u\subset 
\bigcap_{1\leq s<\infty}L^s(\ce^{*},\nu)$.
%In paticular, we have 
%$[\ce]_u\subset L^2(\ce^{*},\nu)\subset [\ce]^{*}_u$.
In addition, for each $1\leq s<\infty$,
the inclusion mapping 
$[\ce]_u\hookrightarrow L^s(\nu)$ is continuous.
\end{proposition}
%\begin{pf}
%Let $\vf\in [\ce]_u$.  Since  
%$[\ce]_u$ is closed under conjugation and 
%pointwise products, $|\vf|^{2q}\in [\ce]_u$
%for any integer $q\geq 1$.  Hence 
%\begin{equation*}
%	\int_{\ce^{*}}|\vf(x)|^{2q}\nu(dx)
%	=\lla \vF_{\nu}, |\vf|^{2q} \rra<\infty.
%\end{equation*}
%Thus $\vf\in L^{2q}(\ce^{*},\nu)$ for any integer $q\geq 1$.
%Therefore it is clear to see that 
%$\vf\in L^p(\ce^{*},\nu)$ for any integer $p\geq 1$.
%\end{pf}
%\bigskip
%\noindent
%{\bf Acknowledgements.} The author is grateful for financial 
%supports from the Daiko Foundation and the Kamiyama Foundation.
%In addition, he wants to give his deepest appreciation
%to Professors I. Kubo, H.-H. Kuo 
%and Dr. U.C. Ji for valuable discussions 
%and constant encouragement.

\begin{example}\label{example-grey}
(Grey noise measure)\\
Let $0<\lambda\leq 1$.
The grey noise measure on $\ce^{*}$ 
is the measure $\nu_{\lambda}$ having
the characteristic function
\begin{equation*}
	L_{\lambda}(|\x|^2_0)
	=\int_{\ce^{*}}e^{i\la x,\x\ra}\nu_{\lambda}(dx),
	\quad \x\in \ce,
\end{equation*}
where $L_{\lambda}(t)$ is the Mittag-Leffler
function with parameter $\lambda$;
\begin{equation*}
	L_{\lambda}(t)=\sum_{n=0}^{\infty}
	{(-t)^n\over\Gamma(1+\lambda n)}.
\end{equation*}
Here $\Gamma$ is the Gamma function.
This measure was introduced by Schneider \cite{schneider}.
It is shown in \cite{kuo96} that
$\nu_{\lambda}$ is a Hida measure which 
induces a generalized function 
$\vF_{\nu_{\lambda}}$ in $(\ce)^{*}_{1-\lambda}$.
Therefore by Theorem \ref{thm:mainthm} and 
Example \ref{example2.8}
the grey noise measure $\nu_{\lambda}$ 
satisfies 
\begin{equation*}
	\int_{\ce^{*}_{p}}
	\exp\bigl({1\over 2}(2-\lambda)
	|x|_{-p}^{{2\over 2-\lambda}}\bigr)\nu_{\lambda}(dx)
	<\infty
\end{equation*}
for some $p$. 
\end{example}

\begin{example}\label{example-poisson}
(Poisson measure)\\
Let ${\cal P}$ be the Poisson measure on $\ce^{*}$ 
given by 
\begin{equation*}
	\exp\Bigl(\int_{\real}(e^{i\x(t)}-1)dt\Bigr)
	=\int_{\ce^{*}}e^{i\la x,\x\ra}{\cal P}(dx),
	\quad \x\in\ce^{*}.
\end{equation*}
It has been shown \cite{cks} that the Poisson noise 
measure induces a generalized function in the Bell 
number space of order $2$.
Thus by Theorem \ref{thm:mainthm} and Example \ref{example2.9}
we have the 
integrability condition
\begin{equation*}
	\int_{\ce^{*}_p}
	\exp\Bigl(|x|_{-p}\sqrt{\log |x|_{-p}}\Bigr)
	{\cal P}(dx)<\infty
\end{equation*}
for some $p$. 
\end{example}

\bigskip
\noindent
{\bf Acknowledgements} 
The author is grateful for financial supports
from the Daiko Foundation and Kamiyama Foundation.
\bigskip
%%%%%%%%%%%%%%%%%%%%%%%%%%%%%%%%%%%%%%%%%%%%%
% References
%%%%%%%%%%%%%%%%%%%%%%%%%%%%%%%%%%%%%%%%%%%%

\end{document}